\documentclass[10pt]{article}

\usepackage{amssymb,amsmath}
\usepackage{color}

\definecolor{pdflinkcolor}{rgb}{.1,.1,.6}	
\definecolor{pdfcitecolor}{rgb}{.6,.1,.1}	

\usepackage[%
	 colorlinks=true, 
	 citecolor=pdfcitecolor,
	 linkcolor=pdflinkcolor,
	 a4paper,
	 pagebackref=true,
	 pdfauthor={B. Jacob, C. Trunk, M. Winklmeier},
	 pdftitle={Analyticity and Riesz basis property of semigroups associated to damped vibrations}]{hyperref}

\newtheorem{thm}{Theorem}[section]
\newtheorem{lem}[thm]{Lemma}

\newtheorem{prop}[thm]{Proposition}
\newtheorem{cor}[thm]{Corollary}
\newtheorem{defn}[thm]{Definition}

\newtheorem{rem}[thm]{Remark}

\makeatletter
\renewcommand*\refstepcounter[1]{\stepcounter{#1}%
  \protected@edef\@currentlabel{%
    \csname p@#1\expandafter\endcsname
      \csname the#1\endcsname
  }}
\makeatother%

\makeatletter
\newcommand{\condition}[1][{{\upshape(\alph{condnumber})}}]{
   \renewcommand{\p@condnumber}[1]{#1}	
   \refstepcounter{condnumber}		
   \item[#1]
}
\makeatother

\newcounter{condnumber}
\newenvironment{conditions}
{
\begin{list}{}{
      \usecounter{condnumber}
      \settowidth{\labelwidth}{(B) }
      \setlength{\labelsep}{1ex}
      \setlength{\leftmargin}{\labelwidth}
      \addtolength{\leftmargin}{\labelsep}
      \setlength{\itemindent}{0pt}
      \setlength{\listparindent}{0pt}
      \setlength{\rightmargin}{0pt}
      \addtolength{\itemsep}{3pt}
   } 
} 
{\end{list}}	

\newcommand{\R}{{\mathbb{R}}}
\newcommand{\N}{{\mathbb{N}}}
\newcommand{\C}{{\mathbb{C}}}
 \newcommand{\Deh}{{\bf {\mathcal D}}}
 \newcommand{\Skindef}{[\raisebox{0.5 ex}{.},\raisebox{0.5 ex}{.}]}
 
 \newcommand{\Norm}{\|\,\raisebox{0.5 ex}{.}\,\|}

 \newcommand{\ov}{\overline}
 \newcommand{\lln}{\lambda_{0}}

\def\cA{{\mathcal A}}      
      
   \def\cH{{\mathcal H}}

\newcommand{\Real}{\mbox{{\rm Re}}\,}

\newcommand{\ds}{\displaystyle}
\newenvironment{mat}{ \left[ \begin{array}{ccccc} }{\end{array} \right]  }

\newcommand{\dom}{{\cal D}}

\newcommand{\la}{\langle}
\newcommand{\ra}{\rangle}

\newcommand{\half}{\frac{1}{2}}

\newcommand{\be}{\begin{equation}}
\newcommand{\ee}{\end{equation}}
\newcommand{\bea}{\begin{eqnarray}}
\newcommand{\eea}{\end{eqnarray}}
\newcommand{\beann}{\begin{eqnarray*}}
\newcommand{\eeann}{\end{eqnarray*}}
\newcommand{\kasten}{{\hspace*{\fill}$\square$}\vspace{1ex}}

\usepackage{fancyhdr}  

\date{}
\title{Analyticity and Riesz basis property of semigroups associated to damped vibrations
}

\author{Birgit Jacob, Carsten Trunk and Monika Winklmeier}

\begin{document}
\maketitle

\begin{abstract} 
   Second order equations of the form $\ddot{z}(t) + A_0 z(t) + D \dot{z} (t)=0$ are considered. 
   Such equations are often used as a model for transverse motions of thin beams in the presence of damping. 
   We derive various properties of the operator matrix $\cA = \left[\begin{smallmatrix} 0 & I \\ -A_0 & -D \end{smallmatrix}\right]$ associated with the second order problem above.
   We develop sufficient conditions for analyticity of the associated semigroup and for the existence of a Riesz basis consisting
of eigenvectors and associated vectors of $\cA$ in the phase space. 
\end{abstract}


\section{Introduction}  
A linear equation describing transverse motions of a thin beam can be written in the form
\begin{equation*}
   \frac{\partial ^2 u } {\partial t^2 } + \frac{\partial
   ^2}{\partial r^2 } \left [ E  \frac{\partial ^2 u } {\partial r^2
   } + {C_d } \frac{\partial ^3 u }{\partial r^2 \partial t } \right] =0,
   \hspace{2em} r \in (0,1), t > 0,
\end{equation*}
where
$u(r,t) $ is the transverse displacement of the beam at time $t$ and position $r$. 
The existence and behaviour of solutions $u$ depend also on boundary and initial conditions.
In the example above we are interested in solution having finite energy, i.e. solutions such that $\|u(\cdot, t)\|^2 + \|u''(\cdot, t)\|^2 < \infty$ for all $t>0$ where $\|\cdot\|$ denotes the usual norm in the Hilbert space $L^2(0,1)$.
Identifying the function $u(\cdot,t)$ with an element $z(t)\in L^2(0,1)$ by $z(t)(r) = u(r,t)$ we obtain from the partial differential equation above a second order equation in $L^2(0,1)$ of the form 
\begin{equation}\ddot{z}(t) + A_0 z(t) + D \dot{z} (t)    =   0,
\label{sys}
\end{equation}
where $A_0 = E \frac{\partial^4}{\partial r^4}$, $D = \frac{\partial^2}{\partial r^2} C_d \frac{\partial^2}{\partial r^2}$ acting in $L^2(0,1)$ with appropriate domains encoding the boundary conditions under consideration.
We will come back to this example in Section~6.

\medskip
In this paper we study second order equations of type \eqref{sys} in an abstract Hilbert space $H$ where the stiffness operator $A_0$ is a possibly unbounded positive operator on $H$ and is assumed to be boundedly invertible, and $D$, the damping operator, is an unbounded operator such that $A_0^{-1/2}D A_0^{-1/2}$ is a bounded non-negative operator on $H$.
This second order equation is equivalent to the standard first-order equation $\dot{x}(t)=\cA x(t)$ where $\cA:\dom(\cA)\subset \dom(A_0^{1/2})\times H\rightarrow \dom(A_0^{1/2})\times H$ 
is given by
$$ \cA = \begin{mat} 0 & I \\ -A_0 & -D \end{mat} , \hspace{1em}$$
$$ \dom(\cA) = \left\{ \left[\begin{smallmatrix} z \\ w \end{smallmatrix}\right] \in
  \dom(A_0^{1/2}) \times \dom(A_0^{1/2}) \mid A_0 z + D w \in H  \right\} . $$
This operator matrix has been studied in the literature for more than 20 years. Interest in this particular model is motivated by various problems such as stabilization, see for example \cite{ben}, \cite{latr}, \cite{lev}, \cite{sle}, solvability of Riccati equations \cite{hela}, minimum-phase property \cite{JMT} and compensator problems with partial observations  \cite{hela2}.

It is well-known that $\cA$ generates a $C_0$-semigroup of contractions
in $H_{1/2}\times H$, where $H_{1/2}=\dom(A_0^{1/2})$ is equipped with the norm $\|x\|_{1/2} = \|A_0^{1/2}x\|_H$, and thus the spectrum of $\cA$ is located in the closed left half plane. 
This goes back to \cite{bi} and \cite{las}, see also \cite{biw}, \cite{CLL}. 
Several authors have proved independently of each other that the condition $  \langle A_0^{- 1/2}Dz, A_0^{1/2}z\rangle_{H} \ge \beta \|z\|^2_H$
is sufficient for exponential stability of the $C_0$-semigroup generated by $\cA$, see for example \cite{bi}, \cite{biw}, \cite{baen}, \cite{CLL},   \cite{HS}, \cite{HSII}, \cite{ves}, \cite{TWII}. 

In this paper we focus on two properties of the operator $\cA$: Analyticity of the generated semigroup and the Riesz basis property in the phase space $H_{1/2}\times H$.
Analyticity of the $C_0$-semigroup generated by $\cA$ has been studied in many papers, see \cite{chru}, \cite{CLL}, \cite{chtr}, \cite{chtr2}, \cite{hrsh}, \cite{hu2}, \cite{hu}, and \cite{lanshk}.
Most of the papers require that the damping operator $D$ is comparable with $A^\rho$ for some $\rho\in[1/2,1]$.  
In \cite{lanshk} the damping operator $D$ is of the form
\begin{equation} \label{a+b}
   D=\alpha A_0 +B,
\end{equation}
where $\alpha>0$ is a constant, $A_0^{-1}$ is compact and $B$ is symmetric and $A_0$-compact. 
If $-1/\alpha\not\in\sigma_p(\cA)$, then it is shown in \cite{lanshk} 
that $\cA$ generates an analytic semigroup. 
In this case the essential spectrum of the  operator $A_0^{-1} D$, considered as an operator acting in $H_{1/2}$ consists of the point $-1/\alpha$ only. 
We extend the result of~\cite{lanshk} to more general damping operators $D$:
If $A_0^{-1}$ is compact  in $H$ and  $0\not\in\sigma_{ess}(A_0^{-1}D)$, then $A$ generates an analytic semigroup  on $H_{1/2}\times H$ (cf.\ Theorem \ref{analytic} below).
In particular, this implies that the above mentioned result from \cite{lanshk} holds even if $-1/\alpha\in\sigma_p(\cA)$.
 
Note that analyticity of the semigroup generated by $\cA$ already implies that the semigroup satisfies the spectral mapping theorem 
\[ \sigma(T(t))\backslash\{0\} = e^{t\sigma(\cA)},\qquad t\ge 0,\]
see \cite[Chapter IV, Section 3.10]{EN}. 

We further develop conditions guaranteeing that the space $H_{1/2} \times H$ possesses
a Riesz basis consisting of eigenvectors and  finitely many associated vectors of $\cA$. 
The existence of such a system has many important implications for the operator $\cA$; for instance, it implies that $\cA$ satisfies the weak spectral mapping theorem, that is, 
\[ \sigma(T(t))=\overline{e^{t\sigma(\cA)}},\quad t\ge 0,\]
where $(T(t))_{t\ge 0}$ is the semigroup generated by $\cA$.
In particular, it follows that the semigroup is exponentially stable if and only if the spectrum of $\cA$ is contained in the open left half plane and uniformly bounded away from the imaginary axis. 

The Riesz basis property has been shown in \cite{lanshk} in the situation where
$A_0^{-1}$ is a compact operator, 
$D$ is of the form (\ref{a+b}) for some $\alpha \geq 0$ with a symmetric operator $B$ and $-1/\alpha\not\in\sigma_p(\cA)$, if $\alpha\neq 0$ 
(and with some additional assumptions in the case $\alpha =0$).
Similar results were obtained \cite[Appendix A]{chtr} 
in a more special situation.
All these assumptions guarantee
 that the essential spectrum of $\cA$ consists at most of one point. 

In this paper we also assume that $A_0^{-1}$ is a compact operator, but we allow a more general damping operators $D$. 
In particular, the essential spectrum of $\cA$ may contain infinitely many points. 
For most of our results we need the assumption that $0\not\in\sigma_{ess}(A_0^{-1}D)$, where $A_0^{-1}D$ is seen as an operator acting on $H_{1/2}=\dom(A_0^{1/2})$. 
This, however, implies that we cannot handle the case where $A_0^{-1}D$ is a compact operator in $H_{1/2}$ unless $H$ is of finite dimension.
Together with some rather weak conditions the above mentioned imply that there exists a Riesz basis in the phase space $H_{1/2}\times H$ consisting of eigenvalues and finitely many associated vectors of $\cA$ (cf.\ Theorem \ref{Riesz} below). 
For results involving compact $A_0^{-1}D$ we refer the reader to \cite{lanshk}.\\

Throughout this paper we assume that all Hilbert and Krein spaces are infinite dimensional.
Since we are interested in applications to partial differential equations, this is no major restriction.\\

We proceed as follows. 
In Section 2 we provide some useful results on the spectrum of operators in Krein spaces. 
In particular, we recall the notion of spectral points of positive and negative type and of type $\pi_+$ and type $\pi_-$. 
One main tool of this paper is to show that certain spectral points of $\cA$ are of positive or negative type or of type $\pi_+$ or $\pi_-$. 
In Section 3 we give the precise definition of the operator $\cA$ and prove some of its properties. 
The main results of this paper are contained in Sections 4 and 5 where we always assume that that $A_0^{-1}$ is a compact operator and that $0\not\in\sigma_{ess}(A_0^{-1}D)$. 
The main result of Section 4 is that $\cA$ generates an analytic strongly continuous semigroup. 
Further, we show that $\infty$ is a spectral point of negative type and that every real spectral point is of type $\pi_+$. 
As a consequence we obtain that $\cA$ is definitizable and that the non-real spectrum of $\cA$ consists of at most finitely many points belonging to the point spectrum of $\cA$. 
Further, the operator $\cA$ can be written as a direct sum of a self-adjoint operator in a Hilbert space and a bounded self-adjoint operator in a Pontryagin space. 
Section 5 is devoted to the Riesz basis property of the operator $\cA$, that is, it is shown that under additional weak conditions there exists a Riesz basis of $\dom(A_0^{1/2}) \times H$ consisting of eigenvalues and finitely many associated vectors of $\cA$. 
Finally, in Section 6 the results are illustrated by an example: the Euler-Bernoulli beam with distributed Kelvin-Voigt damping.



\section{Spectrum of operators in Krein spaces}\label{section1} 

Let $({\mathcal H},\Skindef)$ be a Krein space.  
We briefly recall that a complex linear space $\cH$ with a hermitian nondegenerate sesquilinear form $\Skindef$ is called a {\em Krein space} if there exists a decomposition $\cH = \cH_+ \oplus \cH_-$ with subspaces $\cH_\pm$ being orthogonal to each other with respect to $\Skindef$ such that $(\cH_\pm, \pm\Skindef)$ are Hilbert spaces.
In the following all topological notions are understood with respect to some Hilbert space norm $\Norm$ on ${\cal H}$ such that $\Skindef$ is $\Norm$-continuous.
Any two such norms are equivalent.
For the basic theory of Krein space and operators acting therein we refer to \cite{B} and \cite{AI}.

Let $A$ be a closed operator in ${\mathcal H}$. 
Analogously to the Hilbert space case we define the extended spectrum $\sigma_{e}(A)$ of $A$ by $\sigma_{e}(A):= \sigma(A)$ if $A$ is bounded and 
$\sigma_{e}(A):=\sigma(A) \cup \{\infty \}$ if $A$ is unbounded.
The resolvent set of $A$ is denoted by $\rho(A)$. 
By $\sigma_{p,norm}(A)$ we denote the set of all $\lambda \in {\mathbb C}$ which are isolated points of $\sigma(A)$ and normal eigenvalues of $A$, that is, the corresponding Riesz-Dunford projection is of finite rank.
A point $\lambda_{0} \in \C$ is said to belong to the {\em approximative point spectrum\/} $\sigma_{ap}(A)$ of $A$ if there exists a sequence $(x_{n}) \subset \Deh(A)$  with $\|x_{n} \| =1$, $n= 1,2, \ldots,$ and $\|(A - \lln I)x_{n} \| \to 0$ if $n \to \infty$. 
For a self-adjoint operator $A$ in ${\mathcal H}$ all real spectral points of $A$ belong to $\sigma_{ap}(A)$ (see e.g.\ \cite[Corollary VI.6.2]{B}).
The operator $A$ is called {\em Fredholm} if the dimension of the kernel of $A$ and the codimension of the range of $A$ are finite.
The set
\begin{displaymath}
   \sigma_{ess}(A) := \{ \lambda \in {\mathbb C}\; |\; A-\lambda I\mbox{ is not Fredholm} \}
\end{displaymath}
is called the {\em essential spectrum} of $A$.

Using the indefiniteness of the scalar product on $\cH$ we have the notions of spectral points of positive and negative type and of type $\pi_{+}$ and type $\pi_{-}$.
The following definition was given in \cite{LMM1} and \cite{LMaM} for bounded self-adjoint operators.

\begin{defn}  \label{definition++}
For a  self-adjoint operator $A$ in
${\mathcal H}$
a point $\lln \in \sigma(A)$ is called a
spectral point of {\em positive\/} {\rm (}{\em negative\/}{\rm )}
 {\em type of $A$\/}
if $\lln \in \sigma_{ap}(A)$ and for every sequence 
$(x_{n}) \subset \Deh(A)$ with
$\|x_{n} \| =1$ and $\| (A -\lln I)x_{n} \| \to 0$ as
$n \to \infty$ we have
\begin{displaymath}
\liminf_{n \to \infty}\, [x_{n},x_{n}] >0 \;\;\;\; 
 \bigl(\,resp.\,\, \limsup_{n \to \infty}\, [x_{n},x_{n}] <0 \bigr).
\end{displaymath}
The point $\infty$ is said to be 
 of {\em positive\/} {\rm (}{\em negative\/}{\rm )}
 {\em type of $A$\/} if $A$ is unbounded and
 for every sequence 
$(x_{n}) \subset \Deh(A)$ with
$\lim_{n \to \infty} \|x_{n} \| =0$ and $\| Ax_{n} \| = 1$ we have
\begin{displaymath}
\liminf_{n \to \infty}\, [Ax_{n},Ax_{n}] >0 \;\;\;\;
 \bigl(\,resp.\,\, \limsup_{n \to \infty}\, [Ax_{n},Ax_{n}] <0 \bigr).
\end{displaymath}
We denote the set of all points of $\sigma_{e}(A)$ of positive
{\rm (}negative{\rm )} type by $\sigma_{++}(A)$
 {\rm (}resp.\ $\sigma_{--}(A)${\rm )}. 
\end{defn}

It is not difficult to see that the sets $\sigma_{++}(A)$ 
and $\sigma_{--}(A)$ are contained in $\ov{\R}$.
Moreover the non-real spectrum of $A$ cannot accumulate to 
 $\sigma_{++}(A) \cup \sigma_{--}(A)$.

In a similar way as above we define
subsets $\sigma_{\pi_{+}}(A)$ and $\sigma_{\pi_{-}}(A)$ of $\sigma_{e}(A)$
containing $\sigma_{++}(A)$ and
$\sigma_{--}(A)$, respectively (cf.\ \cite[Definition 5]{AJT}).

\begin{defn} \label{def}
For a  self-adjoint operator $A$ in
${\mathcal H}$  
a point $\lln \in \sigma(A)$ is called a
spectral point of 
{\em type $\pi_{+}$\/} {\rm (}{\em type $\pi_{-}$}\/{\rm )} {\em of $A$\/}
if $\lln \in \sigma_{ap}(A)$ 
and if there exists a linear submanifold ${\mathcal H}_{0} \subset {\mathcal H}$ with
{\rm codim}$\,{\mathcal H}_{0} < \infty$ such that
for every sequence 
$(x_{n}) \subset {\mathcal H}_{0} \cap \Deh(A)$ with
$\|x_{n} \| =1$ and $\| (A -\lln I)x_{n} \| \to 0$ as
$n \to \infty$ we have
\begin{equation}\label{..0}
\liminf_{n \to \infty}\, [x_{n},x_{n}] >0 \;\;\;\; 
 \bigl(\,resp.\,\, \limsup_{n \to \infty}\, [x_{n},x_{n}] <0 \bigr).
\end{equation}
The point $\infty$ is said to be 
of {\em type $\pi_{+}$\/} {\rm (}{\em type $\pi_{-}$\/}{\rm )}
 {\em of $A$\/} if $A$ is unbounded and if
there exists a linear submanifold ${\mathcal H}_{0} \subset {\mathcal H}$ with
{\rm codim}$\,{\mathcal H}_{0} < \infty$ such that
 for every sequence 
$(x_{n}) \subset {\mathcal H}_{0} \cap \Deh(A)$ with
$\lim_{n \to \infty} \|x_{n} \| =0$ and $\| Ax_{n} \| = 1$ we have
\begin{displaymath}
\liminf_{n \to \infty}\, [Ax_{n},Ax_{n}] >0 \;\;\;\; 
 \bigl(\,resp.\,\, \limsup_{n \to \infty}\, [Ax_{n},Ax_{n}] <0 \bigr).
\end{displaymath}
We denote the set of all points of $\sigma_{e}(A)$ of type $\pi_{+}$
{\rm (}type $\pi_{-}${\rm )} of $A$ by $\sigma_{\pi_{+}}(A)$ 
{\rm (}resp.\ $\sigma_{\pi_{-}}(A)${\rm )}. 
\end{defn}
Recall that a self-adjoint operator $A$ in a Krein space 
$({\cal H}, \Skindef)$ is called
{\it definitizable} if
$\rho(A) \ne \emptyset$ and if there exists a 
rational function $p \ne 0$ having poles only in $\rho(A)$
such that $[p(A)x,x] \geq 0$ for all $x \in {\cal H}$.
Then the spectrum of $A$ is real or its non-real part consists of
 a finite number of points. Moreover, $A$ has
a spectral function $E(\raisebox{0.5 ex}{.})$ defined
on the ring generated by all connected subsets of $\ov{\R}$
whose endpoints do not belong to some
finite set which is contained in $\{ t \in \R : p(t) = 0 \}
\cup \{\infty \}$ (see \cite{L}).

For a definitizable operator $A$ a
point $t \in \overline{\R}$ is called a {\it critical point}
of $A$ if there is no
open subset $\Delta$ with $t \in \Delta$ such that either 
$\Delta \subset \sigma_{++}(A)$ or
$\Delta \subset \sigma_{--}(A)$.
A critical point $t$ is called {\it regular} if there exists
an open deleted neighbourhood $\delta_{0}$ of $t$ such that
the set of the projections $E(\delta)$  is bounded where $\delta$ runs
through all intervals $\delta$ with $\overline{\delta} \subset
\delta_{0}$, see \cite{L}.

\begin{thm}\label{Deff}
Let $A$ be a self-adjoint operator in
${\mathcal H}$ satisfying
$$
\sigma_{ess}(A) \subset \R, \quad \infty \in \sigma_{++}(A) \cup \sigma_{--}(A)
\quad \mbox{and} \quad \sigma(A) \cap \R \subset \sigma_{\pi_{+}}(A) \cup \sigma_{\pi_{-}}(A). 
$$
Then $A$ is a definitizable operator, and the non-real spectrum of $\cA$ 
consists of at most finitely many points which belong to
$\sigma_{p,norm}(\cA)$.
\end{thm}
{\bf Proof:}\\
Assume $\infty \in \sigma_{--}(A)$.
By \cite[Lemma 2]{AJT} there exists a neighbourhood ${\cal U}$ of $\infty$
in $\ov{\C}$ with
$$
{\cal U} \setminus \ov{\R} \subset \rho(A) \quad \mbox {and} \quad 
{\cal U} \cap \R \subset \sigma_{--}(A)  \cup \rho(A).
$$
{}From this and \cite[Theorem 18]{AJT} we conclude that the non-real spectrum of $A$
consists of at most finitely many points which belong to 
$\sigma_{p,norm}(A)$. Then, by 
\cite[Theorem 23]{AJT} and \cite[Theorem 4.7]{J2}, the operator $A$
is definitizable. A similar reasoning applies to the case 
$\infty \in \sigma_{++}(A)$.
\kasten


\section{Framework and preliminary results}\label{section2}   

Throughout this paper we make the following assumptions.

{\bf (A1)} The stiffness operator $A_0 : \dom (A_0) \subset H
\rightarrow H$ is a self-adjoint uniformly positive operator
on a Hilbert space $H$. 
We define $H_{\frac{1}{2}} =  \dom(A_0^{1/2}) $ equipped with the norm 
$\| \,\cdot\, \|_{H_{\frac{1}{2}}} := \|A_0^{1/2} \, \cdot\, \|_{H}$
and $H_{-{\frac{1}{2}}}=H_{\frac{1}{2}}^*$. 
Here the duality is taken with respect to the pivot space $H$, that is, equivalently $H_{-{\frac{1}{2}}}$ is the completion of $H$ with respect to the norm $\|z\|_{H_{-{\frac{1}{2}}}}= \|A_0^{-{1/2}}z\|_H$. 
Thus $A_0$ restricts to a bounded operator $A_0:H_{\frac{1}{2}} \rightarrow H_{-{\frac{1}{2}}}$. 
We use the same notation $A_0$ to denote this restriction.

We denote the inner product on $H$ by
$\langle\cdot,\cdot\rangle_H$ or $\langle\cdot,\cdot\rangle$, 
and the duality pairing on $H_{-{\frac{1}{2}}}\times H_{{\frac{1}{2}}}$ by
$\langle\cdot,\cdot\rangle_{H_{-{\frac{1}{2}}}\times H_{{\frac{1}{2}}}}$. 
Note that for $(z',z)\in H\times H_{\frac{1}{2}}$, we have
\[ 
\langle z',z\rangle_{H_{-{\frac{1}{2}}}\times H_{{\frac{1}{2}}}}=\langle z',z\rangle_H.
\]

{\bf (A2)} The  damping operator $D: H_{\frac{1}{2}} \rightarrow
H_{-\frac{1}{2}}$ is a bounded operator such  that
 $A_0^{-1/2}D A_0^{-1/2}$ is a bounded self-adjoint operator in $H$ and satisfies
  \[
 \langle Dz, z\rangle_{H_{-\frac{1}{2}}\times
  H_{\frac{1}{2}}} \ge 0  ,\qquad z\in H_{\frac{1}{2}}.
\]

The equation (\ref{sys})  is equivalent to the
following standard first-order equation
\bea \dot{x} (t) & = & \cA x(t) \label{first1}\eea 
where $\cA:\dom(\cA)\subset
H_{\frac{1}{2}}  \times H\rightarrow H_{\frac{1}{2}}  \times H$,
is given by
$$ \cA = \begin{mat} 0 & I \\ -A_0 & -D \end{mat} , \hspace{1em}$$
$$ \dom(\cA) = \left\{ \left[\begin{smallmatrix} z \\ w \end{smallmatrix}\right] \in
  H_{\frac{1}{2}} \times
H_{\frac{1}{2}} \mid A_0 z + D w \in H  \right\} . $$
The operator $\cA$ itself is not self-adjoint in the Hilbert space $H_{\frac{1}{2}}\times H$. 
It is easy to see (e.g.\ \cite{TWII}) that $\cA$ has a bounded inverse in $H_{\frac{1}{2}}\times H$ given by
\begin{equation} \label{NaDann}
 \cA^{-1}=\begin{mat} -A_0^{-1}D & -A_0^{-1} \\ I & 0 \end{mat},
\end{equation}
where $A_0^{-1}D$ is considered as an operator acting in $H_{\frac{1}{2}}$.
This together with the fact that
\begin{equation*}
 J\cA, \qquad \mbox{ where } J= \begin{mat} I & 0\\ 0 & -I
\end{mat},
\end{equation*}
is a symmetric operator in the Hilbert space $H_{\frac{1}{2}} \times H$, imply the self-adjointness of $J\cA$ in $H_{1/2}\times H$.
Therefore, (compare also \cite[Proof of Lemma 4.5]{WT}) 
\begin{equation*}\label{defj}
 \cA^*= J\cA J,\qquad \mbox{ with } \dom(\cA^*)=J\dom(\cA) 
\end{equation*}
and  
\begin{equation*}
\Real \langle \cA x,x \rangle \leq 0 \quad \mbox{for } x \in \dom(\cA) \quad
\mbox{and} \quad
\Real \langle \cA^* x,x \rangle \leq 0 \quad \mbox{for } x \in \dom(\cA^*).
\end{equation*}
Hence, $\cA$ is the
generator of a strongly continuous semigroup of contractions on
the state space $H_{\frac{1}{2}}  \times H $.
This fact is well-known, see e.g.\ \cite{bi}, \cite{biw},  \cite{CLL}, \cite{HS}, \cite{las} or
\cite{TWII}.

For $\left(\begin{smallmatrix} x_{1}\\
y_{1}\end{smallmatrix}\right),\left(\begin{smallmatrix} x_{2}\\
y_{2}\end{smallmatrix}\right) \in  H_{\frac{1}{2}}\times H$ we define
an indefinite inner product on $H_{\frac{1}{2}}\times H$ by 
\begin{equation} \label{Krein}
\left[ \begin{pmatrix} x_{1}\\
y_{1}\end{pmatrix}, \begin{pmatrix} x_{2}\\
y_{2}\end{pmatrix}\right] := 
\left\la J \begin{pmatrix} x_{1}\\
y_{1}\end{pmatrix}, \begin{pmatrix} x_{2}\\
y_{2}\end{pmatrix}\right\ra = \la x_{1},x_{2} \ra_{H_\frac{1}{2}}-
\la y_{1},y_{2} \ra.
\end{equation}
Then $(H_{\frac{1}{2}}\times H, \Skindef )$ is a Krein space  
and $\cA$ is a self-adjoint operator
 with respect to $\Skindef$.

In the following proposition we collect the above considerations. 
 
\begin{prop}  \label{thm:cont}
The operator $\cA$ is self-adjoint in the Krein space
$(H_{\frac{1}{2}}\times H, \Skindef )$, its spectrum
is contained
in the closed left half plane and lies symmetric 
with respect to the real line. 
The operator $\cA$ 
has a bounded inverse and is the
generator of a strongly continuous semigroup of contractions on
$H_{\frac{1}{2}}  \times H $.
\end{prop}

This implies that the spectrum of $\cA$ is a subset of the closed left half plane without the origin and symmetric with respect to the real axis.
However, otherwise the spectrum of $\cA$ is quite arbitrary. 
For an example with $\sigma(\cA)=\{ s\in \mathbb C\mid {\rm Re}\,s\le 0, |s|\ge \varepsilon\}$, $\varepsilon>0$, 
 we refer to \cite{JT}.

In the following theorem
we give an estimate for the neighbourhood of the origin which lies in the resolvent set
of $\cA$ and for the modulus of
the eigenvalues of $\cA$.

\begin{thm} \label{Winkl}
We have
$\lambda \in \rho(\cA)$ if and only if the operator 
$I + \lambda A_0^{-1}(D + \lambda I)$, considered as an
operator in ${\cal L}(H_{\frac{1}{2}})$, is boundedly invertible. In particular
\begin{equation} \label{Winkl1}
\left\{ \lambda \in \C  \; |\; \|\lambda A_0^{-1}D + 
\lambda^2 A_0^{-1}\|_{{\cal L}(H_{\frac{1}{2}})} < 1 \right\}
\subset \rho(\cA). 
\end{equation}
Moreover, each $\lambda \in \sigma_{p}(\cA)$ satisfies
\begin{align}\label{eq:EVestimate}
   |\lambda | \geq \frac{1}{2\|A_0^{-1}\|_{{\cal L}(H)}}\left( \sqrt{
	 \|A_0^{-1}D\|_{{\cal L}(H_{\frac{1}{2}})}^2 + 4\|A_0^{-1}\|_{{\cal L}(H)}} - 
      \|A_0^{-1}D\|_{{\cal L}(H_{\frac{1}{2}})}
   \right).
\end{align}
\end{thm}
{\bf Proof:}\\
Let $\lambda \in \rho(\cA)$. Then by \cite[Proposition 2.2]{JT}
the operator 
$$
 \lambda^2 A_0^{-1} +\lambda A_0^{-1/2} D A_0^{-1/2} + I
$$
is bounded and  boundedly invertible in $H$, hence
$$
 A_0^{-1/2} \left(\lambda^2 A_0^{-1} +\lambda A_0^{-1/2} D A_0^{-1/2} + I \right) A_0^{1/2}
= I + \lambda A_0^{-1}(D + \lambda I)
$$
is boundedly invertible in $H_{\frac{1}{2}}$. For the contrary choose
$\lambda \in \C$, $\left(\begin{smallmatrix} u\\
v \end{smallmatrix}\right) \in \Deh(\cA)$ and
$\left( \begin{smallmatrix} x\\
y\end{smallmatrix}\right) \subset H_{\frac{1}{2}} \times H$. Then we have
$$
(\cA - \lambda I) \left( \begin{smallmatrix} u\\
v \end{smallmatrix}\right) = 
\left(\begin{smallmatrix} x\\
y\end{smallmatrix}\right) 
$$
if and only if the following equations hold
\begin{eqnarray*}
& v = x + \lambda u &\\
& -A_0 (I + \lambda A_0^{-1}(D + \lambda I))u = y +Dx + \lambda x. &
\end{eqnarray*}
This implies the first assertion  of Theorem \ref{Winkl}.
 Let $\lambda \in \sigma_{p}(\cA)$.
Then the above calculations imply that the operator $I+\lambda A_0^{-1}(D+\lambda I)$
in $H_{\frac{1}{2}}$ is not injective. Therefore, there exists a non-zero
vector $f \in H_{\frac{1}{2}}$ with
$$
f = -\lambda A_0^{-1}(D+\lambda I)f.
$$
Hence,
$$
\|f\|_{H_{\frac{1}{2}}} \leq |\lambda | \left(
 \|A_0^{-1}D\|_{{\cal L}(H_{\frac{1}{2}})} +|\lambda | 
\|A_0^{-1}\|_{{\cal L}(H_{\frac{1}{2}})}
\right)\|f\|_{H_{\frac{1}{2}}}
$$
and, as  $\|A_0^{-1}\|_{{\cal L}(H_{\frac{1}{2}})} = \|A_0^{-1}\|_{{\cal L}(H)}$, we conclude
$$
\left( |\lambda | + \frac{\|A_0^{-1}D\|_{{\cal L}(H_{\frac{1}{2}})}}{ 2\|A_0^{-1}\|_{{\cal L}(H)}}
\right)^2 - \frac{1}{ \|A_0^{-1}\|_{{\cal L}(H)}} - 
\frac{\|A_0^{-1}D\|^2_{{\cal L}(H_{\frac{1}{2}})}}{ 4\|A_0^{-1}\|^2_{{\cal L}(H)}} \geq 0
$$
and Theorem \ref{Winkl} is proved.
\kasten

\begin{rem}
   The estimate~\eqref{eq:EVestimate} for the eigenvalues is optimal since in the case $D=0$ it
    follows that  $\mu$ is an eigenvalue of $A_0$ if and only if $\pm i \sqrt{\mu}$ are
     eigenvalues of $\cA$.
   If the uniformly positive operator $A_0$ has a compact resolvent, then the 
    smallest eigenvalue of $A_0$ equals $\| A_0^{-1} \|^{-1} $ and the eigenvalue
    $\lambda_{\min}$ of $\cA$ with smallest
    absolute eigenvalue is given by
   \begin{align*}
      |\lambda_{\min}| 
      \, =\, 
      \sqrt{ \min\{ \mu\, | \, \mu \text{ eigenvalue of } A_0 \} }
      \, =\, 
      \sqrt{ \| A_0^{-1} \|^{-1} }
   \end{align*}
   which is equal to the right hand side of~\eqref{eq:EVestimate} if $D$ is set to be $0$.

\end{rem}



\section{Analyticity}\label{section3}   

Throughout this section we assume that $A_0^{-1}$ is a compact operator.
Note that $A_0^{-1} D$, considered as an operator acting in $H_{\frac{1}{2}}$, is a bounded  non-negative operator. 
In \cite[Theorem 4.1]{JT} it is shown that under this assumption for $\lambda\in \mathbb C\backslash\{0\}$ we have 
\begin{equation}\label{ess}
   \lambda\in \sigma_{ess}(-A_0^{-1}D)\qquad \mbox{if and only if }
   \quad 1/\lambda\in \sigma_{ess}(\cA).
\end{equation}

If not explicitely stated otherwise, the operator $A_0^{-1}D$ is always considered as an operator acting on $H_{1/2}$.

We obtain the following main result concerning analyticity.

\begin{thm}\label{analytic}
Assume that $A_0^{-1}$ is compact  in $H$ and that $0\not\in\sigma_{ess}(A_0^{-1}D)$. Then $A$ generates an analytic semigroup  on $H_{1/2}\times H$.
\end{thm}

The proof of this theorem will be given at the end of this section. We first prove some properties of the point infinity. The following theorem shows in particular that $\infty \in \sigma_{--}(\cA)$ if
$0\not\in \sigma_{ess}(A_0^{-1}D)$.

\begin{thm} \label{infty}
Assume that the operator $A_0^{-1}$ is a compact operator in $H$ and that $0\not\in \sigma_{ess}(A_0^{-1}D)$. Then 
$$
\infty \in \sigma_{--}(\cA) \quad \mbox{and} \quad \R \subset \sigma_{\pi_{+}}(\cA)
\cup \rho(\cA).
$$
Moreover, 
the operator $\cA$ is definitizable and
there exists a neighbourhood ${\cal U}$ of $\infty$ in
$\ov{\C}$ and constants $M>0$, $m \in \N$ and $\eta>0$  such that
\begin{equation} \label{ess2}
{\cal U} \setminus \ov{\R} \subset \rho(\cA) \quad \mbox {and} \quad 
{\cal U} \cap \R \subset \sigma_{--}(\cA) \cup \rho(\cA)
\end{equation}
and
\begin{eqnarray} \label{ess3}
&\|(\cA-\lambda I)^{-1}\| \leq \frac{M}{|{\rm Im}\, \lambda|} \quad \mbox{for all }
\lambda \in {\cal U}\setminus \ov{\R},& \\ \label{ess4}
&\|(\cA-\lambda I)^{-1}\| \leq \frac{M}{|{\rm Im}\, \lambda|^m} \quad \mbox{for all }
\lambda \in \rho(A) \setminus \R \mbox{ with }|{\rm Im}\lambda|\le \eta.
\end{eqnarray}
Further, the non-real spectrum of $\cA$ 
consists of at most finitely many points which belong to
$\sigma_{p,norm}(\cA)$.
\end{thm}
{\bf Proof:}\\
The proof is divided into two steps. 
First we will prove that $\infty \in \sigma_{--}(\cA)$.
In the second step we will show that and $\R \subset \sigma_{\pi_{+}}(\cA) \cup \rho(A)$.
Since by \eqref{ess} the essential spectrum of $\cA$ is real, Theorem \ref{Deff} yields that $\cA$ is a definitizable operator and the non-real spectrum of $\cA$ 
consists of at most finitely many points which belong to
$\sigma_{p,norm}(\cA)$.
Further, (\ref{ess2}), (\ref{ess3}) and (\ref{ess4}) follow from
\cite[Lemma 2 and Proposition 3]{AJT} and from
\cite[Proposition II.2.1]{L}. 

\medskip
\emph{Step 1.} By \cite[Lemma 10]{AJT}, $\infty$ belongs to
$\sigma_{--}(\cA)$ if and only if $\infty$ belongs to  $\sigma_{\pi_-}(\cA)$. 
It is easily seen (see e.g.\ \cite{ABJT}) that this is the case if and only if  $0 \in \sigma_{\pi_-}(\cA^{-1})$.
Assume $0 \notin \sigma_{\pi_-}(\cA^{-1})$. 
Then there exists a sequence
$\left(\left(\begin{smallmatrix} x_{n}\\
y_{n}\end{smallmatrix}\right)\right) \subset H_{\frac{1}{2}} \times H$ 
with $\left\|\left(\begin{smallmatrix} x_{n}\\
y_{n}\end{smallmatrix}\right)\right\|^2_{H_{\frac{1}{2}}\times H} = \|x_n\|^2_{H_{\frac{1}{2}}} +
\|y_{n}\|^2 =1$ and $\cA^{-1} \left(\begin{smallmatrix} x_{n}\\
y_{n}\end{smallmatrix}\right) \to 0$ as $n \to \infty$ and
\begin{equation} \label{grross}
   \limsup_{n \to \infty} \left[ 
      \left(\begin{smallmatrix} x_{n}\\ y_{n}\end{smallmatrix}\right),
      \left(\begin{smallmatrix} x_{n}\\ y_{n}\end{smallmatrix}\right)
   \right] 
   = \limsup_{n \to \infty} 
   \left( \|x_n\|^2_{H_{\frac{1}{2}}} -  \|y_n\|^2 \right)
   \geq 0.
\end{equation}
By \cite[Theorem 14]{AJT} this sequence can be chosen to converge to zero weakly.
This gives
\begin{equation} \label{grross2}
\|A_0^{-1}Dx_{n} + A_0^{-1}y_n\|_{H_{\frac{1}{2}}} \to 0 \quad \mbox{and} \quad
\|x_{n}\| \to 0 \quad \mbox{as } n \to \infty.
\end{equation}
The sequence $(A_0^{-1/2}y_{n})$ converges weakly to zero in $H_{\frac{1}{2}}$.
As $A_0^{-1}$ is a compact operator in $H$, $A_0^{-1/2}$ is a compact operator in $H_{\frac{1}{2}}$. 
It follows that $(A_0^{-1}y_{n})$ converges to zero in $H_{\frac{1}{2}}$.
Then, by (\ref{grross2}), we have
$$
\|A_0^{-1}Dx_{n}\|_{H_{\frac{1}{2}}} \to 0  \quad \mbox{as } n \to \infty.
$$
Moreover, the sequence $(x_{n})$ converges weakly to zero in
$H_{\frac{1}{2}}$, hence the assumption $0 \notin \sigma_{ess}(A_0^{-1}D)$ implies
$$
\|x_{n}\|_{H_{\frac{1}{2}}} \to 0  \quad \mbox{as } n \to \infty.
$$
Then $\|y_n\| \to 1$ as $n \to \infty$, in contradiction to (\ref{grross}) and $0 \in \sigma_{\pi_-}(\cA^{-1})$ follows.

\medskip
\emph{Step 2.}
We now choose $\mu \in (-\infty,0)$ and 
$$
G_{\mu} := \mbox{span}\, \{ x \in H_{\frac{1}{2}} \mid A_0 x = \nu x, \nu \leq 
\mu^2 \}.
$$
Then $G_{\mu}$ is a finite dimensional subspace of $H_{\frac{1}{2}}$.
For every sequence $\left(\left(\begin{smallmatrix} x_n\\
y_n\end{smallmatrix}\right)\right) $ in $\dom(\cA) \cap (G_{\mu}
\times G_{\mu})^{\perp}$ with
$\left\|\left(\begin{smallmatrix} x_{n}\\
y_{n}\end{smallmatrix}\right)\right\|^2_{H_{\frac{1}{2}}\times H}=1$ and
 $(\cA-\mu I) \left(\begin{smallmatrix} x_{n}\\
y_{n}\end{smallmatrix}\right) \to 0$ as $n \to \infty$ we have
$$
\|y_n - \mu x_n \|_{H_{\frac{1}{2}}} \to 0 \quad \mbox{and} \quad
\|A_0 x_n + Dy_n + \mu y_n\| \to 0  \quad \mbox{as } n \to \infty.
$$
This gives
\begin{equation*}
\begin{split}
\liminf_{n \to \infty}
\left[ \left(\begin{smallmatrix} x_{n}\\
y_{n}\end{smallmatrix}\right),\left(\begin{smallmatrix} x_{n}\\
y_{n}\end{smallmatrix}\right)\right]& = \liminf_{n \to \infty}
\left( \la x_{n},x_{n} \ra_{H_\frac{1}{2}} -\la y_{n},y_{n} \ra\right)\\
& = \liminf_{n \to \infty} \left(\la A_0x_{n},x_{n} \ra
- \mu^{2} \la x_{n},x_{n} \ra\right) >0,
\end{split}
\end{equation*}
where the last inequality follows from the fact that $x_n \in G_{\mu}^\perp$,
$n \in \N$. Therefore $\R \subset \sigma_{\pi_{+}}(\cA)$ and Theorem 
\ref{infty} is proved.
\kasten

\begin{rem}
The stronger assumption $0\not\in\sigma(A_{0}^{-1}D)$ implies 
that there exist constants $\alpha, \gamma >0$ with
$$
\gamma
\la A_{0}x,x \ra_{H_{-\frac{1}{2}} \times H_{\frac{1}{2}}} \leq 
\la D x,x \ra_{H_{-\frac{1}{2}} \times H_{\frac{1}{2}}} \leq 
\alpha \la A_{0}x,x \ra_{H_{-\frac{1}{2}} \times H_{\frac{1}{2}}}
\quad \mbox{for } x \in H_{\frac{1}{2}}.
$$
\end{rem}

\noindent
{\bf Proof of Theorem \ref{analytic}:} \\
Since $\cA$ is the generator of a strongly continuous semigroup, estimate~\eqref{ess3} shows immediately that $\cA$ generates an analytic semigroup, see~\cite[Chapter II, Section 4.5]{EN}.
\kasten

The following corollary shows that under the assumptions of Theorem \ref{infty} the operator $\cA$ can be written as a
direct sum of a self-adjoint operator on a Hilbert space and a bounded self-adjoint operator on a Pontryagin space. In the situation of $D=\rho A_0^\alpha$, $\rho>0$ and $\alpha\in(0,1]$, $\cA$ is the direct sum of two normal operators \cite{chtr}.

\begin{cor}
Assume that the operator $A_0^{-1}$ is a compact operator in $H$ and that $0\not\in \sigma_{ess}(A_0^{-1}D)$.
Then the  space $({H_{\frac{1}{2}}}\times H, \Skindef)$ decomposes into the direct sum of two $\cA$-invariant
closed subspaces $H^{\prime}$ and $H^{\prime \prime}$, which are orthogonal with respect to $\Skindef$, such that:
\begin{enumerate}
\item The space $(H^{\prime}, -\Skindef)$ is a Hilbert space, $\cA |H^{\prime}$ is a self-adjoint
operator in this Hilbert space and
$$
\sigma(\cA |H^{\prime}) \subset \ov{\R} \setminus (-M,\infty),
$$
where $M$ is as in $(${\rm \ref{umg}}$)$.
\item The space $(H^{\prime \prime}, \Skindef)$ is a Pontryagin space, $\cA |H^{\prime \prime}$ is a bounded self-adjoint
operator in  this Pontryagin space with 
$$
\sigma(\cA |H^{\prime \prime}) \subset [-M,0) \cup \Theta \quad \mbox{and} \quad
\sigma(\cA |H^{\prime \prime}) \subset \sigma_{++}(\cA |H^{\prime \prime}) \cup \Xi
$$
where $\Xi, \Theta \subset \C$ are  empty or consist of finitely many
points and  $\Theta \subset \sigma_{p,norm}(\cA |H^{\prime \prime})$.
\end{enumerate}
\end{cor}
{\bf Proof:}\\
By Theorem \ref{infty} the operator $\cA$ is definitizable with
$\infty \in \sigma_{--}(\cA)$ and $\R \subset \sigma_{\pi_{+}}(\cA) \cup \rho(\cA)$.
Denote by $E$ the spectral function of $\cA$.
Since $\infty \in \sigma_{--}(\cA)$, there exists $M>0$ with 
\begin{equation} \label{umg}
\ov{\R} \setminus
(-M,M) \subset \sigma_{--}(\cA) \cup \rho(\cA),
\end{equation}
see \cite[Lemma 2]{AJT}. 
Set $\Delta_0 := \ov{\R} \setminus [-M,M]$, $H^{\prime}:=E(\Delta_0)
 (H_{\frac{1}{2}}\times H)$ and $H^{\prime \prime}:=(I-E(\Delta_0))
 (H_{\frac{1}{2}}\times H)$. Then the assertions above follow from \cite[Theorem 3.18]{J2}
and Theorem \ref{thm:cont}.
\kasten

\begin{rem}
Note that the essential spectrum of $\cA$ is empty if and only if either the essential spectrum 
of $A_0^{-1}D$ is zero or empty, see \eqref{ess}.
\end{rem}


\section{Expansion in eigenfunctions}\label{section4}   

In the sequel we always assume the Hilbert space $H$ to be separable.
An at most countably infinite set ${\cal M}$ of elements
of a Hilbert space is said to be a {\em Riesz basis} if there exists 
an isomorphic mapping ${\cal M}$ onto an orthonormal basis, cf.\
\cite[Lecture VI]{N}. 

Condition (\ref{strong}) of Theorem \ref{Riesz} below appears
already in the celebrated works \cite{KL1,KL2}, where the case
of a bounded self-adjoint operator $D$ and a positive compact operator
$A_0$ was discussed. We will use this approach in the proof of Theorem \ref{Riesz}.

\begin{thm}\label{Riesz}
Assume that the operator $A_0^{-1}$ is compact in $H$ and that 
\begin{equation}\label{Nulll}
0\not\in \sigma_{ess}(A_0^{-1}D),
\end{equation}
where $A_0^{-1}D$ is considered as an operator acting in
$H_{\frac{1}{2}}$. 
Assume that the set $\sigma_{ess}(A_0^{-1}D)$ is countably and has
at most countable many accumulation points.
Moreover, let at least one of the following conditions be satisfied.

\begin{conditions}
   
   \condition \label{item:condition1}
   There exists a $\delta >0$ such that for all
   $f \in H_{\frac{1}{2}}$ with $\|f\|_{H_{\frac{1}{2}}}=1$ we have
   \begin{equation}\label{strong}
      \left \la A_0^{-1}Df,f \right \ra^2_{H_{\frac{1}{2}}} -
      4 \left \la A_0^{-1}f,f \right \ra_{H_{\frac{1}{2}}} > \delta.
   \end{equation}

   \condition \label{item:condition2}
   For all $\mu \in \sigma_{ess}(-A_0^{-1}D)$ we have either
   $\frac{1}{\mu} \notin \sigma_p(\cA)$ or, if $\frac{1}{\mu} \in \sigma_p(\cA)$,
   there exists no non-zero 
   $\left(\begin{smallmatrix} y\\
	 \mu^{-1} y \end{smallmatrix}\right)
   \in$ {\rm ker}$\, (\cA-\mu^{-1} I)$ such that
   \begin{equation} \label{nondeg}
      \mu^2 \la y,w\ra_{H_{\frac{1}{2}}} = \la y,w\ra 
      \quad \mbox{for all }
      \left(\begin{smallmatrix} w\\
	    \mu^{-1} w \end{smallmatrix}\right)
      \in \mbox{\rm ker} \, (\cA-\mu^{-1}I).
   \end{equation}

   \condition \label{item:condition3}
   $\|A_0^{-1/2}\| <\inf \{\lambda>0\mid \lambda \in \sigma_{ess}(A_0^{-1}D)\}$.
\end{conditions}

Then the following assertions hold.

\begin{enumerate}

   \item \label{item:1}
   There exists a subspace of $H_{\frac{1}{2}} \times H$ of at most finite codimension
    which has a Riesz basis consisting of eigenvectors of $\cA$.

   \item \label{item:2}
   There exists a Riesz basis of $H_{\frac{1}{2}} \times H$ consisting
   of eigenvectors and finitely many associated vectors of $\cA$. 
   
   \item \label{item:4}
   Moreover, if \ref{item:condition1} holds, then $\cA$ has no associated vectors, i.e.\ there are no
   Jordan chains of length greater than one, the spectrum of $\cA$ is real 
   and there exists a Riesz basis 
   of $H_{\frac{1}{2}} \times H$ consisting
   of eigenvectors of $\cA$.
\end{enumerate}
\end{thm}

As a corollary we obtain the following result

\begin{cor} \label{cor52}
Assume that the operator $D$ has the form
\[ D=\alpha A_0 +B,\]
where $\alpha>0$ is a constant and $B$ is a symmetric $A_0$-compact operator.
 If $-1/\alpha\not\in\sigma_p(\cA)$, then there exists a Riesz basis of $H_{\frac{1}{2}} \times H$ consisting
of eigenvectors and finitely many associated vectors of $\cA$, and  $\cA$ generates an analytic semigroup.
\end{cor}
{\bf Proof:}\\
We define $H_1 =  \dom(A_0) $ equipped with the norm $\| \,\cdot\, \|_{H_1} := \|A_0 \, \cdot\, \|$ and $H_{-1}$ is the completion of $H$ with respect to the norm $\| \,\cdot\, \|_{H_{-1}} := \|A_0^{-1} \, \cdot\, \| $. 
Then, by assumption, $B$, restricted from $H_1$ to $H$, is a compact  
operator, cf.\ \cite[IV 1.12]{K}. 
By the symmetry of $B$, $B^*$ is an extension of $B$ and a compact operator
acting from $H$ into $H_{-1}$. Thus, by interpolation,
the operator $B$ considered as an operator from $H_{\frac{1}{2}}$ to $H_{-\frac{1}{2}}$ is compact,
hence $\sigma_{ess}(A_0^{-1}D) = \alpha$ and Corollary \ref{cor52} follows from Theorem \ref{Riesz}.
\kasten

In \cite{lanshk} it is shown that under the assumption of the corollary there exists a Riesz basis of $H_{\frac{1}{2}} \times H_{\frac{1}{2}}$ consisting
of eigenvectors and finitely many associated vectors of $\cA$, and  $\cA$ generates an analytic semigroup.
\vspace{2ex}

\noindent{\bf Proof of Theorem \ref{Riesz}:}\\
It suffices to prove Part~\ref{item:1} and \ref{item:4} of the theorem, as Part~\ref{item:2} follows immediately from Part~\ref{item:1}.

\medskip
Assume that condition \ref{item:condition1} holds.
By Theorem \ref{infty} we have $\R \subset \sigma_{\pi_{+}}(\cA)
\cup \rho(A)$.
Let $\lambda \in\sigma_{\pi_{+}}(\cA)$. Then there exists a sequence 
$\left(\left(\begin{smallmatrix} x_n\\
y_n\end{smallmatrix}\right)\right) $ in $\dom(\cA)$ with
$\left\|\left(\begin{smallmatrix} x_{n}\\
y_{n}\end{smallmatrix}\right)\right\|^2_{H_{\frac{1}{2}}\times H}=1$ and
 $(\cA-\lambda I) \left(\begin{smallmatrix} x_{n}\\
y_{n}\end{smallmatrix}\right) \to 0$ as $n \to \infty$. 
Hence
$$
\|y_n - \lambda x_n \|_{H_{\frac{1}{2}}} \to 0 \quad \mbox{and} \quad
\|A_0x_n + \lambda Dx_n + \lambda^2 x_n\| \to 0 
 \quad \mbox{as } n \to \infty.
$$
This implies
\begin{equation} \label{NNeu}
\|x_n + \lambda A_0^{-1}Dx_n + \lambda^2 A_0^{-1} x_n\|_{H_{\frac{1}{2}}} \to 0 
 \quad \mbox{as } n \to \infty.
\end{equation}
and
\begin{equation*}
\begin{split}
\liminf_{n \to \infty}
\left[ \left(\begin{smallmatrix} x_{n}\\
y_{n}\end{smallmatrix}\right),\left(\begin{smallmatrix} x_{n}\\
y_{n}\end{smallmatrix}\right)\right]& = \liminf_{n \to \infty}
\left( \la x_{n},x_{n} \ra_{H_\frac{1}{2}} -\la y_{n},y_{n} \ra\right)\\
& = \liminf_{n \to \infty} \left(\la x_{n},x_{n} \ra_{H_\frac{1}{2}}
- \lambda^{2} \la x_{n},x_{n} \ra \right)\\
& = \liminf_{n \to \infty} -\lambda \left(\la A_0^{-1}Dx_{n},x_{n} \ra_{H_{\frac{1}{2}}}
+ 2 \lambda \la A_0^{-1}x_{n},x_{n} \ra_{H_{\frac{1}{2}}}\right).
\end{split}
\end{equation*}
Similarly, we have
\begin{equation*}
\limsup_{n \to \infty}
\left[ \left(\begin{smallmatrix} x_{n}\\
y_{n}\end{smallmatrix}\right),\left(\begin{smallmatrix} x_{n}\\
y_{n}\end{smallmatrix}\right)\right]
= \limsup_{n \to \infty} -\lambda \left(\la A_0^{-1}Dx_{n},x_{n} \ra_{H_{\frac{1}{2}}}
+ 2 \lambda \la A_0^{-1}x_{n},x_{n} \ra_{H_{\frac{1}{2}}}\right) .
\end{equation*}
Now \ref{item:condition1} implies $\sigma(\cA) \subset \R$ (see, e.g.\ \cite[Theorem 3.3]{JT}),
hence, by Proposition~\ref{thm:cont}, we have $\lambda \in (-\infty,0)$. 
Moreover, by \ref{item:condition1}, the operator pencil
$$
L(s) := s^2I+sA_0^{-1}D + A_0^{-1}, \quad s \in \C,
$$
considered as a pencil with values in the bounded operators acting on
$H_{\frac{1}{2}}$, is strongly hyperbolic, see e.g.\ \cite[Lemma 31.23]{M}.
Therefore, see e.g.\ \cite{L68}, we have
$$\liminf_{n \to \infty}
\left( \la A_0^{-1}Dx_{n},x_{n} \ra_{H_{\frac{1}{2}}}
+ 2 \lambda \la A_0^{-1}x_{n},x_{n} \ra_{H_{\frac{1}{2}}}\right) >0$$
or
$$\limsup_{n \to \infty}
\left( \la A_0^{-1}Dx_{n},x_{n} \ra_{H_{\frac{1}{2}}}
+ 2 \lambda \la A_0^{-1}x_{n},x_{n} \ra_{H_{\frac{1}{2}}}\right)<0.$$
This gives
$$
\sigma(\cA) \subset \sigma_{++}(\cA) \cup \sigma_{--}(\cA) \subset \R
$$
and $\cA$ has no associated vectors 
and there exists a Riesz basis 
of $H_{\frac{1}{2}} \times H$ consisting
of eigenvectors of $\cA$.

\medskip
Assume that condition \ref{item:condition2} holds and let $\mu \in \sigma_{ess}(-A_0^{-1}D)$ such that $\frac{1}{\mu} \in \sigma_p(\cA)$.
Now, \ref{item:condition2} implies that there are no Jordan chains of $\cA$ corresponding
to the eigenvalue $\frac{1}{\mu}$ of length greater than one, and, moreover, that
ker$\, (\cA - \mu^{-1} I)$ is a non-degenerate subspace of 
$(H_{\frac{1}{2}} \times H, \Skindef)$, that is, 
$$
\mbox{ker}\, (\cA - \mu^{-1} I)  \cap 
\left( \mbox{ker}\, (\cA - \mu^{-1} I) \right)^{[\perp ]}=\{0\},
$$
where $\left( \mbox{ker}\, (\cA - \mu^{-1} I) \right)^{[\perp ]}$ is the 
orthogonal companion of $\mbox{ker}\, (\cA - \mu^{-1} I)$ with respect to
$\Skindef$. Moreover,  (\ref{nondeg}) implies that $\mu^{-1}$
is a regular critical point of $\cA$, see \cite{L} or \cite[Proposition 1.4]{DL}.
 As all points from $\sigma(\cA) \setminus 
\sigma_{ess}(\cA)$ belong to $\sigma_{p,norm}(\cA)$, it turns out that
$\cA$ has only regular critical points and the eigenvectors
of $\cA$ form a Riesz basis of a subspace of $H_{\frac{1}{2}} \times H$
of an at most finite codimension.  The eigenvectors and associated vectors
of $\cA$ form a Riesz basis of $H_{\frac{1}{2}} \times H$. 

\medskip 
Condition \ref{item:condition3} implies \ref{item:condition2}, hence
Theorem~\ref{Riesz} is proved. \kasten

Theorem \ref{Riesz} implies that the operator $\cA$ is the direct sum of an operator similar to a self-adjoint operator in a Hilbert space and a bounded operator in a finite-dimensional space.

\begin{rem}
Assume that the operator $A_0^{-1}$ is compact in $H$ and that {\rm (\ref{Nulll})} holds. 
Then it was shown in the proof of Theorem {\rm \ref{Riesz}} that
if {\rm (i)} holds or if for all $\mu^{-1} \in \sigma_p(\cA)$
we have that {\rm (\ref{nondeg})} holds, all critical points of $\cA$ are regular and there are no associated vectors.
Hence, $\cA$ is similar to a self-adjoint operator in the Hilbert space $H_{\frac{1}{2}} \times H$. 
\end{rem}

\begin{rem}
   We mention that Theorem {\rm \ref{Riesz}} can be obtained also by methods from {\rm \cite{AI}}. 
   For this, one has to show that the operator $\cA^{-1}$ belongs to the class $({\bf H})$, cf.\ {\rm \cite[Chapter 3, \S 5]{AI}}, and then apply {\rm \cite[Theorem 4.2.12]{AI}}.
\end{rem}


\section{Example: Euler-Bernoulli Beam with distri\-bu\-ted Kelvin-Voigt damping}   

We consider a beam of length 1 with a thin film of piezoelectric polymer applied to one side and
we study transverse vibrations only. 
Let $u(r,t) $ denote the deflection of the beam  from its rigid body motion at time $t$ and position $r$. 
Use of the Euler-Bernoulli model for the beam deflection and the Kelvin-Voigt damping model leads to the following description of the vibrations \cite{BSY}, \cite{Pu}:
\begin{equation}
   \frac{\partial ^2 u } {\partial t^2 } + \frac{\partial
   ^2}{\partial r^2 } \left [ E  \frac{\partial ^2 u } {\partial r^2
   } + {C_d } \frac{\partial ^3 u }{\partial r^2 \partial t } \right] =0,
   \hspace{2em} r \in (0,1), t > 0.
   \label{beam}
\end{equation}//
Here the flexural rigidity $E$ is a positive physical constant which is determined by the beam's area momentum of inertia and its modulus of elasticity and $C_d\in L^\infty(0,1)$ with $C_d(t)\ge c>0$ a.e.~describes the damping properties of the piezoelectric film.
Assuming that the beam is pinned at point 0 and sliding at 1, we have for all $t>0$ the  following boundary conditions:
\begin{equation}
u\big|_{r=0} =\   0,  \quad
 \frac{\partial u }{\partial r}\bigg|_{r=1}  =\ 0,  \quad
   \frac{\partial ^2 u } {\partial
      r^2 }\bigg|_{r=0}  =\  0,\quad 
  \frac{\partial ^3 u } {\partial
      r^3 }\bigg|_{r=1}  =\  0.
   \label{bcs}
\end{equation} 
We consider the partial differential equation \eqref{beam}-\eqref{bcs} as a second order problem in the Hilbert space $H = L^2 (0,1)$. In $H$ we define the operator $A_0$ by  
\begin{align*}
   A_0 = E \frac{d^4}{dr^4},
   \quad
   \dom(A_0) = \left\{ z \in H^4(0,1) \mid z(0)=
   z'(1)=  z'' (0)= z''' (1) =0 \right\}. 
\end{align*}
It is easy to see that the operator $A_0$ satisfies  assumption (A1) and that $A_0^{-1}$ is a compact operator. 
We have
$$
H_{\half} = \left\{ z \in H^2(0,1)\mid z(0)=
z'(1) =0 \right\}
$$
with inner product $  \la z, v \ra_{H_\half} = E\la z'', v'' \ra $. The operator $A_0^{1/2}$ is given by 
\begin{equation} \label{Delft}
 A_0^{1/2} = E^{1/2}\frac{d^2}{dr^2} \quad \mbox{and} \quad \|z\|_{H_{\frac{1}{2}}}^2 \geq \frac{\pi^4E}{16} \|z\|^2
\quad \mbox{for } z \in  H_{\frac{1}{2}}.
\end{equation}
Let $\ds x(t) = (u(\cdot,t) , \dot{u} (\cdot,t))$.
Then $\| x(t) \|_{H_{\half}\times H}^2 = \| u''(\cdot, t)\|^2 +  \| \dot{u} (\cdot, t)\|^2 $ corresponds to the energy of the beam which justifies the choice of $L^2(0,1)$ as the Hilbert space for the analysis of the boundary value problem \eqref{beam}-\eqref{bcs}.

By $M_{C_d}\in {\cal L}(H)$ we denote the multiplication operator
\[ (M_{C_d}f)(x) = C_d(x)f(x)\]
and we define the damping operator as
\begin{align*} 
   D=\frac{1}{E} A_0^{1/2}M_{C_d}A_0^{1/2}.
\end{align*}
For $z\in H_\half$ we have
\begin{equation}\label{eqn1} 
   \langle Dz,z\rangle_{H_{-\half}\times H_{\half}} 
   = \langle C_dz'',z''\rangle_H 
   \ge c \|z''\|^2 \ge \frac{\pi^4}{16} c \|z\|^2,
\end{equation}
and thus the assumption (A2) holds as well. 
Furthermore, each solution of the abstract problem $ \ddot{z}(t) + A_0 z(t) + D \dot{z} (t) = 0$ corresponds to a solution of the boundary value problem \eqref{beam}-\eqref{bcs}.
We have the following lemma.

\begin{lem}\label{lemess}
The operators $A_0^{-1/2}M_{C_d}A_0^{-1/2}$ and $A_0^{-1}D$ are bounded self-adjoint operators in $H$ and $H_\half$, respectively, with
\[ \sigma_{ess}(A_0^{-1}D)=\sigma_{ess}(A_0^{-1/2}M_{C_d}A_0^{-1/2})=\sigma(E^{-1} M_{C_d}).\]
\end{lem}
{\bf Proof:}
A vector $x\in H_\half$ belongs to ker$(A_0^{-1}D-\lambda I)$ if and only if  $A_0^{1/2} x$ belongs to ker$(A_0^{-1/2}DA_0^{-1/2}-\lambda I)$. The operator $A_0^{1/2}$ maps $H_\half$ isometrically onto $H$, therefore
\[ {\rm dim}\, {\rm ker}\,(A_0^{-1}D-\lambda I)= {\rm dim}\, {\rm ker}\,(A_0^{-1/2}DA_0^{-1/2}-\lambda I).\]
Obviously, $A_0^{-1/2} D A_0^{-1/2} = M_{C_d}$ and $A_0^{-1}D$ are bounded self-adjoint operators in $H$ and $H_\half$, respectively, and therefore we have for real $\lambda$
\begin{eqnarray*}
 {\rm codim}\, {\rm ran}\,(A_0^{-1}D-\lambda I) &=& {\rm dim}\, {\rm ker}\,(A_0^{-1}D-\lambda I) \\
 &=& {\rm dim}\, {\rm ker}\,(A_0^{-1/2}DA_0^{-1/2}-\lambda I)\\
 &=& {\rm codim}\, {\rm ran}\,(A_0^{-1/2}DA_0^{-1/2}-\lambda I).
\end{eqnarray*}
\vspace{-7.5ex}\\
\kasten

\smallskip
Lemma \ref{lemess} implies $0\not\in\sigma_{ess}(A_0^{-1}D)$, as the function $C_d$ satisfies $C_d(t)\ge c>0$ a.e.
By (\ref{eqn1}), the corresponding operator $\cA$ generates an exponentially stable
semigroup on $ H_{\half} \times L^2(0,1)$ (see the introduction). Moreover, the assumptions of Theorem \ref{analytic}  are satisfied and thus $\cA$ generates an analytic semigroup. By Theorem \ref{infty}, $\cA$ is definitizable, $\infty \in \sigma_{--}(\cA)$, $\R \subset \sigma_{\pi_{+}}(\cA)\cup \rho(A)$ and the non-real spectrum of $\cA$ consists of at most finitely many points which belong to $\sigma_{p,norm}(\cA)$.

In addition, we now assume that the film on the beam consists of several patches, that is,
\[ C_d(x)=\sum_{k=1}^n a_k \chi_{A_k}(x),\]
where $n\in\mathbb N$, $a_k>0$, $k=1,\cdots,n$, $A_k$ are measurable disjoint subsets of $(0,1)$ and 
\[ \overline{\bigcup_{k=1}^n A_k} = [0,1].\]

\begin{thm} \label{endee}
If
\begin{equation}\label{akk}
a_k > \frac{8}{\pi^2\sqrt{E}} 
\end{equation}
holds for $k=1,\ldots,n$, then \ref{item:condition1}  from Theorem  \ref{Riesz} is satisfied, that is,
 the spectrum of $\cA$ is real  and there exists a Riesz basis 
of $H_{\frac{1}{2}} \times H$ consisting of eigenvectors of $\cA$.
If for all $a_k$, $k=1,\ldots,n$, with
$$
a_k \leq \frac{4}{\pi^2}\sqrt{E} 
$$
we have $-E/a_k \notin \sigma_p(\cA)$, then  \ref{item:condition2} of Theorem  \ref{Riesz} is satisfied, that is,
there exists a Riesz basis 
of $H_{\frac{1}{2}} \times H$ consisting of of eigenvectors and finitely many associated vectors of $\cA$.
In particular, this holds true if we have
$$
a_k > \frac{4}{\pi^2}\sqrt{E}  \quad \mbox{for } k=1,\ldots,n.
$$
\end{thm}
{\bf Proof:}\\
By (\ref{Delft}) and (\ref{akk}), we have for $f \in H_{\frac{1}{2}}$
\begin{align*}
\langle A_0^{-1}Df,f\rangle^2_{H_{\frac{1}{2}}} &= \langle Df,f\rangle^2_{H_{-\frac{1}{2}} \times H_{\frac{1}{2}}}
\geq \left( \min \{a_k \mid k=1,\ldots,n \} \right)^2 \|f\|^4_{H_{\frac{1}{2}}} \\
&>4 \|f\|^2 \|f\|^2_{H_{\frac{1}{2}}} + \delta \|f\|^4_{H_{\frac{1}{2}}}
\end{align*}
for some sufficiently small $\delta >0$ and the first assertion of Theorem \ref{endee}
is proved.

The second and third assertion follow from the fact that for $a_k > \frac{4}{\pi^2}\sqrt{E} $
we have 
$$
a_k \langle y,y\rangle_{H_{\frac{1}{2}}} > \|y\|^2
$$
and \ref{item:condition2} (resp.\ \ref{item:condition3}) of Theorem  \ref{Riesz} is satisfied.
\kasten

There is an obvious generalization of Theorem \ref{endee} for the case of countably many patches which we do not give here in detail.




\begin{thebibliography}{10}

\def\cprime{$'$} \def\cprime{$'$}

\bibitem{ABJT}
T.Ya. Azizov, J.~Behrndt, P.~Jonas, and C.~Trunk.
\newblock Spectral points of type $\pi_+$ and $\pi_-$ for closed linear
  relations in {K}rein spaces.
\newblock submitted.

\bibitem{AI}
T.Ya. Azizov and I.S. Iokhvidov.
\newblock {\em Linear operators in spaces with an indefinite metric}.
\newblock Pure and Applied Mathematics (New York). John Wiley \& Sons Ltd.,
  Chichester, 1989.

\bibitem{AJT}
T.Ya. Azizov, P.~Jonas, and C.~Trunk.
\newblock Spectral points of type {$\pi\sb +$} and {$\pi\sb -$} of self-adjoint
  operators in {K}rein spaces.
\newblock {\em J. Funct. Anal.}, 226(1):114--137, 2005.

\bibitem{bi}
H.T. Banks and K.~Ito.
\newblock A unified framework for approximation in inverse problems for
  distributed parameter systems.
\newblock {\em Control Theory Adv. Tech.}, 4(1):73--90, 1988.

\bibitem{biw}
H.T. Banks, K.~Ito, and Y.~Wang.
\newblock Well posedness for damped second-order systems with unbounded input
  operators.
\newblock {\em Differential Integral Equations}, 8(3):587--606, 1995.

\bibitem{BSY}
H.T. Banks, R.C. Smith, and Y.~Wang.
\newblock The modeling of piezoceramic patch interactions with shells, plates,
  and beams.
\newblock {\em Quart. Appl. Math.}, 53(2):353--381, 1995.

\bibitem{baen}
A.~B{\'a}tkai and K.~Engel.
\newblock Exponential decay of {$2\times2$} operator matrix semigroups.
\newblock {\em J. Comput. Anal. Appl.}, 6(2):153--163, 2004.

\bibitem{ben}
C.D. Benchimol.
\newblock A note on weak stabilizability of contraction semigroups.
\newblock {\em SIAM J. Control Optimization}, 16(3):373--379, 1978.

\bibitem{B}
J.~Bogn{\'a}r.
\newblock {\em Indefinite inner product spaces}.
\newblock Springer-Verlag, New York, 1974.
\newblock Ergebnisse der Mathematik und ihrer Grenzgebiete, Band 78.

\bibitem{chru}
G.~Chen and D.L. Russell.
\newblock {A mathematical model for linear elastic systems with structural
  damping.}
\newblock {\em Q. Appl. Math.}, 39:433--454, 1982.

\bibitem{CLL}
S.~Chen, K.~Liu, and Z.~Liu.
\newblock Spectrum and stability for elastic systems with global or local
  {K}elvin-{V}oigt damping.
\newblock {\em SIAM J. Appl. Math.}, 59(2):651--668 (electronic), 1999.

\bibitem{chtr}
S.~Chen and R.~Triggiani.
\newblock Proof of extensions of two conjectures on structural damping for
  elastic systems.
\newblock {\em Pacific J. Math.}, 136(1):15--55, 1989.

\bibitem{chtr2}
S.~Chen and R.~Triggiani.
\newblock {Characterization of domains of fractional powers of certain
  operators arising in elastic systems, and applications.}
\newblock {\em J. Differ. Equations}, 88(2):279--293, 1990.

\bibitem{DL}
A.~Dijksma and H.~Langer.
\newblock Operator theory and ordinary differential operators.
\newblock In {\em Lectures on operator theory and its applications (Waterloo,
  ON, 1994)}, volume~3 of {\em Fields Inst. Monogr.}, pages 73--139. Amer.
  Math. Soc., Providence, RI, 1996.

\bibitem{EN}
K.~Engel and R.~Nagel.
\newblock {\em One-parameter semigroups for linear evolution equations}, volume
  194 of {\em Graduate Texts in Mathematics}.
\newblock Springer-Verlag, New York, 2000.

\bibitem{hrsh}
R.O. Griniv and A.A. Shkalikov.
\newblock {Operator models in elasticity theory and hydromechanics and the
  associated analytic semigroups.}
\newblock {\em Mosc. Univ. Math. Bull.}, 54(5):1--10, 1999.

\bibitem{HS}
R.O. Griniv and A.A. Shkalikov.
\newblock Exponential stability of semigroups related to operator models in
  mechanics.
\newblock {\em Mat. Zametki}, 73(5):657--664, 2003.

\bibitem{hela}
E.~Hendrickson and I.~Lasiecka.
\newblock Numerical approximations and regularizations of {R}iccati equations
  arising in hyperbolic dynamics with unbounded control operators.
\newblock {\em Comput. Optim. Appl.}, 2(4):343--390, 1993.

\bibitem{hela2}
E.~Hendrickson and I.~Lasiecka.
\newblock Finite-dimensional approximations of boundary control problems
  arising in partially observed hyperbolic systems.
\newblock {\em Dynam. Contin. Discrete Impuls. Systems}, 1(1):101--142, 1995.

\bibitem{HSII}
R.O. Hryniv and A.A. Shkalikov.
\newblock {Exponential decay of solution energy for equations associated with
  some operator models of mechanics.}
\newblock {\em Funct. Anal. Appl.}, 38(3):163--172, 2004.

\bibitem{hu2}
F.~Huang.
\newblock {On the mathematical model for linear elastic systems with analytic
  damping.}
\newblock {\em SIAM J. Control Optimization}, 26(3):714--724, 1988.

\bibitem{hu}
F.~Huang.
\newblock {Some problems for linear elastic systems with damping.}
\newblock {\em Acta Math. Sci.}, 10(3):319--326, 1990.

\bibitem{JMT}
B.~Jacob, K.~Morris, and C.~Trunk.
\newblock Minimum-phase infinite-dimensional second-order systems.
\newblock To appear in {\em IEEE Transactions on Automatic Control}.

\bibitem{JT}
B.~Jacob and C.~Trunk.
\newblock Location of the spectrum of operator matrices which are associated to
  second order equations.
\newblock {\em Operators and Matrices}, 1:45--60, 2007.

\bibitem{J2}
P.~Jonas.
\newblock On locally definite operators in {K}rein spaces.
\newblock In {\em Spectral analysis and its applications}, volume~2 of {\em
  Theta Ser. Adv. Math.}, pages 95--127. Theta, Bucharest, 2003.

\bibitem{K}
T.~Kato.
\newblock {\em Perturbation theory for linear operators}.
\newblock Springer-Verlag, Berlin, second edition, 1976.
\newblock Grundlehren der Mathematischen Wissenschaften, Band 132.

\bibitem{KL1}
M.G. Krein and H.~Langer.
\newblock {On some mathematical principles in the linear theory of damped
  oscillations of continua I.}
\newblock {\em Integral Equations Oper. Theory}, 1:364--399, 1978.

\bibitem{KL2}
M.G. Krein and H.~Langer.
\newblock {On some mathematical principles in the linear theory of damped
  oscillations of continua II.}
\newblock {\em Integral Equations Oper. Theory}, 1:539--566, 1978.

\bibitem{LMM1}
P.~Lancaster, A.S. Markus, and V.I. Matsaev.
\newblock Definitizable operators and quasihyperbolic operator polynomials.
\newblock {\em J. Funct. Anal.}, 131(1):1--28, 1995.

\bibitem{lanshk}
P.~Lancaster and A.A. Shkalikov.
\newblock Damped vibrations of beams and related spectral problems.
\newblock {\em Canad. Appl. Math. Quart.}, 2(1):45--90, 1994.

\bibitem{L68}
H.~Langer.
\newblock \"{U}ber stark ged\"ampfte {S}charen im {H}ilbertraum.
\newblock {\em J. Math. Mech.}, 17:685--705, 1967/1968.

\bibitem{L}
H.~Langer.
\newblock Spectral functions of definitizable operators in {K}re\u\i n spaces.
\newblock In {\em Functional analysis (Dubrovnik, 1981)}, volume 948 of {\em
  Lecture Notes in Math.}, pages 1--46. Springer, Berlin, 1982.

\bibitem{LMaM}
H.~Langer, A.S. Markus, and V.I. Matsaev.
\newblock Locally definite operators in indefinite inner product spaces.
\newblock {\em Math. Ann.}, 308(3):405--424, 1997.

\bibitem{las}
I.~Lasiecka.
\newblock Stabilization of wave and plate-like equations with nonlinear
  dissipation on the boundary.
\newblock {\em J. Differential Equations}, 79(2):340--381, 1989.

\bibitem{latr}
I.~Lasiecka and R.~Triggiani.
\newblock Uniform exponential energy decay of wave equations in a bounded
  region with {$L\sb 2(0,\infty; L\sb 2(\Gamma))$}-feedback control in the
  {D}irichlet boundary conditions.
\newblock {\em J. Differential Equations}, 66(3):340--390, 1987.

\bibitem{lev}
N.~Levan.
\newblock The stabilizability problem: a {H}ilbert space operator decomposition
  approach.
\newblock {\em IEEE Trans. Circuits and Systems}, 25(9):721--727, 1978.
\newblock Special issue on the mathematical foundations of system theory.

\bibitem{M}
A.S. Markus.
\newblock {\em Introduction to the spectral theory of polynomial operator
  pencils}, volume~71 of {\em Translations of Mathematical Monographs}.
\newblock American Mathematical Society, Providence, RI, 1988.

\bibitem{N}
N.K. Nikol{\cprime}ski{\u\i}.
\newblock {\em Treatise on the shift operator}, volume 273 of {\em Grundlehren
  der Mathematischen Wissenschaften [Fundamental Principles of Mathematical
  Sciences]}.
\newblock Springer-Verlag, Berlin, 1986.

\bibitem{sle}
M.~Slemrod.
\newblock Stabilization of boundary control systems.
\newblock {\em J. Differential Equations}, 22(2):402--415, 1976.

\bibitem{WT}
M.~Tucsnak and G.~Weiss.
\newblock How to get a conservative well-posed linear system out of thin air.
  {II}. {C}ontrollability and stability.
\newblock {\em SIAM J. Control Optim.}, 42(3):907--935 (electronic), 2003.

\bibitem{ves}
K.~Veseli{\'c}.
\newblock Energy decay of damped systems.
\newblock {\em ZAMM Z. Angew. Math. Mech.}, 84(12):856--863, 2004.

\bibitem{TWII}
G.~Weiss and M.~Tucsnak.
\newblock How to get a conservative well-posed linear system out of thin air.
  {I}. {W}ell-posedness and energy balance.
\newblock {\em ESAIM Control Optim. Calc. Var.}, 9:247--274 (electronic), 2003.

\bibitem{Pu}
P.H. You.
\newblock Boundary feedback control of elastic beam equation with structural
  damping and stability.
\newblock {\em Acta Math. Appl. Sinica (English Ser.)}, 6(4):373--382, 1990.

\end{thebibliography}

\ \\

\newlength{\adressenbreite}
\settowidth{\adressenbreite}{\sl Sekretariat MA 6-3, Stra\ss e des 17.~Juni 136 }
\hfill\parbox{\adressenbreite}{\sl
   Birgit Jacob\\
   Department of Applied Mathematics\\
   Delft University of Technology\\
   P.O.~Box 5031, 2600 GA Delft\\
   The Netherlands\\
   e-mail: b.jacob@tudelft.nl\\
}

   \medskip
\hfill\parbox{\adressenbreite}{\sl
   Carsten Trunk\\
   Institut f\"ur Mathematik\\
   Technische Universit\"at Berlin\\
   Sekretariat MA 6-3, Stra\ss e des 17.~Juni 136\\
   D-10623 Berlin, Germany\\
   e-mail: trunk@math.tu-berlin.de\\
}

   \medskip
\hfill\parbox{\adressenbreite}{\sl
   Monika Winklmeier\\
   Mathematisches Institut,  Universit\"at Bern\\
   Sidlerstrasse 5\\
   CH-3012 Bern\\
   Switzerland\\
   e-mail: monika.winklmeier@math.unibe.ch
}

\end{document}